\documentclass[11pt]{amsart}
\usepackage{amsmath}
\usepackage{amscd}
\usepackage{graphics}
\usepackage{latexsym}
\oddsidemargin .in
\theoremstyle{plain}
\newtheorem{Thm}{Theorem}
\newtheorem{Lem}[Thm]{Lemma}
\newtheorem{Cor}[Thm]{Corollary}
\newtheorem{Prop}[Thm]{Proposition}

\newtheorem{Def}{Definition}

\theoremstyle{remark}

\def\E{{\mathbb E}}
\def\R{{\mathbb R}}

\def\O{{\mathbb O}}
\def\H{{\mathbb H}}
\def\V{{\mathbb V}}

\def\CL{\mathcal L}

\def\s2x{\hbox{$S^2 \times S^2$}}

    \def\sqr#1#2{{\vcenter{\hrule height.#2pt
            \hbox{\vrule width.#2pt height#1pt \kern#1pt
            \vrule width.#2pt}\hrule height.#2pt}}}
    \def\square{\mathchoice\sqr67\sqr67\sqr{2.1}6\sqr{1.5}6}

\def\qed{~\hfill$\square$}

\begin{document}

\title[]{Harvey Lawson Manifolds and Dualities}
\author{Selman Akbulut and Sema Salur}
\thanks{S. Akbulut is partially supported by NSF grant DMS 0505638, S.Salur is partially supported by NSF grant 1105663}
\address{Department  of Mathematics, Michigan State University, East Lansing, MI, 48824}
\email{akbulut@math.msu.edu }
\address {Department of Mathematics, University of Rochester, Rochester, NY, 14627 }
\email{salur@math.rochester.edu } 
\date{\today}

\begin{abstract}
The purpose of this paper is to introduce Harvey-Lawson manifolds and review the construction of  certain ``mirror dual'' Calabi-Yau submanifolds inside a $G_2$ manifold.  More specifically, given a Harvey-Lawson manifold $HL$, we explain how to assign a pair of tangent bundle valued 2 and 3-forms to a $G_2$
manifold $(M,HL, \varphi, \Lambda)$, with the calibration 3-form
$\varphi$ and an oriented  $2$-plane field $\Lambda$.  As in \cite{as2} these
forms can then be used  to define different complex and
symplectic structures on certain  6-dimensional subbundles of
$T(M)$. When these bundles are integrated they give mirror CY
manifolds (related thru HL manifolds).
\end{abstract}

\maketitle

\setcounter{section}{0}
\vspace{-0.3in}

\section{Introduction}

Let $(M^7,\varphi)$ be a $G_2$ manifold with the calibration 3-form
$\varphi$. If $\varphi$ restricts to be the volume form of an oriented
3-dimensional submanifold $Y^3$, then $Y$ is called an associative
submanifold of $M$.  In \cite{as2} the authors introduced a notion of mirror duality in any $G_2$ manifold 
$(M^{7}, \varphi)$ based on  the associative/coassociative splitting of its tangent bundle $TM=\E\oplus \V$ by the non-vanishing  $2$-fame fields provided by \cite{t}. This duality initially depends  on the choice of two non-vanishing vector fields, one in $\E$ and the other in $\V$.  In this article we give a natural form of this duality where the choice of these vector fields are made more canonical, in the expense of possibly localizing this process to the tubular neighborhood of the $3$-skeleton of $(M,\varphi)$.

\vspace{.07in}

\section{Basic Definitions}

Let us recall some basic facts about $G_2$ manifolds  (e.g. \cite{b1}, \cite{hl}, \cite{as1}). Octonions give an 8 dimensional division algebra $\O=\H\oplus l \H= \R^8$ generated by  $\langle 1, i, j, k, l, li ,lj, lk\rangle $. The imaginary octonions $im \O =\R^7$ is equipped with the cross product operation $\times:  \R^7\times \R^7 \to \R^7$ defined  by $u\times v=im(\bar{v} .u)$. The exceptional Lie group $G_{2}$ is
the linear automorphisms of  $im \O$ preserving this cross product.  Alternatively: 
\begin{equation}
G_{2} =\{ (u_{1},u_{2},u_{3})\in (\R^{7})^3 \;|\; \langle u_{i},u_{j} \rangle=\delta_{ij}, \; \langle
u_{1} \times u_{2},u_{3} \rangle =0 \;  \}.
\end{equation}
\begin{equation} G_{2}=\{ A \in GL(7,\R) \; | \; A^{*} \varphi_{0} =\varphi_{0}\; \}. 
\end{equation}
\vspace{.01in}
\noindent where 
$\varphi _{0}
=e^{123}+e^{145}+e^{167}+e^{246}-e^{257}-e^{347}-e^{356}$ with
$e^{ijk}=dx^{i}\wedge dx^{j} \wedge dx^{k}$.
We say a $7$-manifold $M^7$ has a {\it $G_{2}$ structure} if there is  a
3-form $\varphi \in \Omega^{3}(M)$ such that  at each $p\in  M$
the pair $ (T_{p}(M), \varphi (p) )$ is (pointwise) isomorphic to $(T_{0}(
\R^{7}), \varphi_{0})$. This condition is equivalent to reducing the tangent frame bundle of $M$ from $GL(7, \R)$ to $G_2$.
A manifold with $G_{2}$ structure $(M,\varphi)$  is called a
{\it $G_{2}$  manifold} (integrable $G_2$ structure) if at each point $p\in M$ there is
a chart  $(U,p) \to (\R^{7},0)$ on which $\varphi $ equals to
$\varphi_{0}$ up to second order term, i.e. on the image of the open set $U$ we can write $\varphi (x)=\varphi_{0} + O(|x|^2)$. 
\vspace{.1in}

One important class of $G_2$ manifolds are the ones obtained
from Calabi-Yau manifolds. Let $(X,\omega, \Omega)$ be a complex
3-dimensional Calabi-Yau manifold with K\"{a}hler form $\omega$
and a nowhere vanishing holomorphic 3-form $\Omega$, then
$X^6\times S^1$ has holonomy group $SU(3)\subset G_2$, hence is a
$G_2$ manifold. In this case $ \varphi$= Re $\Omega + \omega
\wedge dt$. Similarly, $X^6\times \mathbb{R}$ gives a noncompact
$G_2$ manifold.

\vspace{.05in}

{\Def Let $(M, \varphi )$ be a $G_2$ manifold. A 4-dimensional
submanifold $X\subset M$ is called {\em coassociative } if
$\varphi|_X=0$. A 3-dimensional submanifold $Y\subset M$ is called
{\em associative} if $\varphi|_Y\equiv vol(Y)$; this condition is
equivalent to the condition $\chi|_Y\equiv 0$,  where $\chi \in \Omega^{3}(M,
TM)$ is the  tangent bundle valued 3-form defined by the
identity:}
\begin{equation}
\langle \chi (u,v,w) , z \rangle=*\varphi  (u,v,w,z)
\end{equation}
The equivalence of these  conditions follows from  the `associator equality' of  \cite{hl}
\begin{equation*}
\varphi  (u,v,w)^2 + |\chi (u,v,w)|^2/4= |u\wedge v\wedge w|^2
\end{equation*}

Similar to the definition of $\chi$ one can define a tangent
bundle 2-form, which is just the cross product of $M$
(nevertheless viewing it as a $2$-form has its advantages).

\vspace{.05in}

{\Def Let $(M, \varphi )$ be a $G_2$ manifold. Then $\psi \in
\Omega^{2}(M, TM)$ is the tangent bundle valued 2-form defined by
the identity:}
\begin{equation}
\langle \psi (u,v) , w \rangle=\varphi  (u,v,w)=\langle u\times v , w
\rangle
\end{equation}

\vspace{.05in}

On a local chart of a $G_2$ manifold $(M, \varphi )$, the form $\varphi $ coincides with the form $\varphi_{0} \in \Omega^{3} (\R^7)$  up to quadratic terms, we can express the tangent valued forms $\chi$ and $\psi$ in terms of $\varphi_{0}$ in local coordinates. More generally, if $e_1,...e_7$ is any local orthonormal frame and $e^1,..., e^7$ is the dual frame, from definitions we get:

 \begin{equation*}
\begin{aligned}
\chi=&\;\;(e^{256}+e^{247}+e^{346}-e^{357})e_1\\
&+(-e^{156}-e^{147}-e^{345}-e^{367})e_2\\
&+(e^{157}-e^{146}+e^{245}+e^{267})e_3\\
&+(e^{127}+e^{136}-e^{235}-e^{567})e_4\\
&+(e^{126}-e^{137}+e^{234}+e^{467})e_5\\
&+(-e^{125}-e^{134}-e^{237}-e^{457})e_6\\
&+(-e^{124}+e^{135}+e^{236}+e^{456})e_7.\\
\end{aligned}
\end{equation*}

\vspace{.1in}

\begin{equation*}
\begin{aligned}
\psi=&\;\;(e^{23}+e^{45}+e^{67})e_1\\
&+(e^{46}-e^{57}-e^{13})e_2\\
&+(e^{12}-e^{47}-e^{56})e_3\\
&+(e^{37}-e^{15}-e^{26})e_4\\
&+(e^{14}+e^{27}+e^{36})e_5\\
&+(e^{24}-e^{17}-e^{35})e_6\\
&+(e^{16}-e^{25}-e^{34})e_7.\\
\end{aligned}
\end{equation*}

Here are some useful facts :
\vspace{.1in}

{\Lem (\cite{as1})  To any $3$-dimensional submanifold $Y^3\subset
(M,\varphi)$,  $\chi$ assigns a normal vector field, which
vanishes when $Y$ is associative.}

{\Lem(\cite{as1}) To any associative manifold $Y^3\subset (M,\varphi)$ with a non-vanishing oriented $2$-plane field,  $\chi$ defines a complex structure on  its normal bundle (notice in particular that any coassociative submanifold $ X\subset M$ has an almost complex structure if its normal bundle has a non-vanishing section).}

\begin{proof}
Let $L\subset \R^7$ be an associative $3$-plane,
that is $\varphi_{0}|_{L}=vol(L)$. Then for every pair of
orthonormal vectors $\{u,v\}\subset L$, the form $\chi$ defines a
complex structure on the orthogonal $4$-plane $L^{\perp}$, as
follows:   Define $j: L^{\perp} \to L^{\perp}$ by
\begin{equation}
j(X)=\chi(u,v,X)
\end{equation}
This is well defined i.e. $j(X)\in L^{\perp}$, because when $ w\in L$ we have:
$$ \langle \chi(u,v,X),w \rangle=*\varphi_{0}(u,v,X,w)=-*\varphi_{0}(u,v,w,X)=\langle \chi(u,v,w),X \rangle=0$$

\noindent Also $j^{2}(X)=j(\chi(u,v,X))=\chi(u,v,\chi(u,v,X))=-X$. We can check the last equality by taking an orthonormal basis $\{ X_{j}\}\subset L^{\perp}$ and calculating
\begin{eqnarray*}
\langle \chi(u,v,\chi(u,v,X_{i})),X_{j}\rangle &=&*\varphi_{0}(u,v,\chi(u,v,X_{i}),X_{j})=-*\varphi_{0}(u,v,X_{j},\chi(u,v,X_{i}))\\
&=&- \langle \chi(u,v,X_{j}),\chi(u,v,X_{i})\rangle =-\delta_{ij}
\end{eqnarray*}

The last equality holds since the map $j$ is orthogonal, and the
orthogonality  can be seen by polarizing the associator equality, and by noticing $\varphi_{0}(u,v,X_i)=0$. Observe that the
map $j$ only depends on the oriented $2$-plane $\Lambda=<u,v>$ generated by $\{u,v\}$ (i.e. it only depends on the complex structure on
$\Lambda $). \end{proof}

\section{Calabi-Yau's hypersurfaces in $G_2$ manifolds}

In \cite{as2}  authors proposed a notion of {\it mirror duality} for Calabi-Yau submanifold pairs lying inside of a $G_{2}$ manifold $(M, \varphi)$. This is done first by assigning almost Calabi-Yau structures to hypersurfeces induced  by hyperplane distributions. The construction goes as follows.  Suppose $\xi$ be a nonvanishing vecor fileld $\xi \in \Omega^{0}(M,TM)$, which gives a codimension one integrable distribution $V_{\xi}:= \xi^{\perp}$ on $M$. If $X_{\xi}$ is a leaf of this distribution, then the forms $\chi$ and $\psi$ induce a non-degenerate $2$-form $\omega_{\xi}$ and an almost complex structure $J_{\xi}$ on $X_{\xi}$ as follows: 
\begin{equation}
\omega_\xi=\langle \psi, \xi \rangle \;\;\;\mbox{and} \;\; J_{\xi}(u)=u\times \xi.
\end{equation}
\begin{equation}
 \textup{Re}\; \Omega_{\xi} = \varphi|_{V_{\xi}}
\;\;\mbox{and}\;\;\; \textup{Im}\; \Omega_{\xi} = \langle \chi, \xi \rangle.
\end{equation}
\noindent where the inner products, of the vector valued differential forms  $\psi$ and $\chi$ with vector field $\xi$, are performed by using their vector part.  So $\omega_{\xi}=\xi \lrcorner \; \varphi$, and $\textup{Im}\;\Omega_{\xi} =\xi \lrcorner \; *\varphi $. Call $\Omega_{\xi} =  \textup{Re}\; \Omega_{\xi} +i \;\textup{Im}\; \Omega_{\xi}$. These induce almost Calabi-Yau structure on  $X_{\xi}$,    
analogous to Example 1.


\begin{Thm} (\cite{as2}) Let $(M,\varphi)$ be a $G_2$ manifold, and $\xi $ be a  unit
vector field such that $\xi^{\perp}$ comes  from a codimension one foliation on $M$,
then  $(X_{\xi},\omega_{\xi}, \Omega_{\xi},J_{\xi})$ is an almost
Calabi-Yau manifold such that $\varphi |_{X_{\xi}}= Re\; \Omega_{\xi} $
and  $*\varphi |_{X_{\xi}}= *_{3}\; \omega_{\xi} $. Furthermore, if
$\CL_{\xi}(\varphi )|_{X_{\xi}}=0$ then $d\omega_{\xi}=0$,  and if
$\CL_{\xi}(*\varphi)|_{X_{\xi}}=0$ then $J_{\xi}$ is integrable; when both conditions are satisfied  $(X_{\xi},\omega_{\xi}, \Omega_{\xi},J_{\xi})$  is a Calabi-Yau
manifold.
\end{Thm}


Here is a brief discussion of \cite{as2} with explanation of its terms: Let  ${\xi}^{\#}$ be the dual $1$-form of $\xi$, and  $e_{\xi^{\#}}$  and $i_{\xi}=\xi \lrcorner $ denote the exterior and interior product operations on differential forms, then 
\begin{equation*}
\varphi =e_{\xi^{\#}}\circ i_{\xi }(\varphi )+i_{\xi }\circ e_{\xi^{\#}}(\varphi )=\omega_{\xi}\wedge  \xi^{\#}+Re \; \Omega_{\xi}.
\end{equation*}
This is the decomposition of the form $\varphi $ with respect to $\xi \oplus \xi^{\perp}$.
The condition that the distribution $V_{\xi}$ to be integrable is $d{\xi}^{\#}\wedge {\xi}^{\#}=0$.
 Also it is clear from definitions that  $J_{\xi}$ is an almost complex
structure on $X_{\xi}$, and the $2$-form $\omega_{\xi }$ is
non-degenerate on $X_{\xi}$, because 
\begin{equation*}
\omega_{\xi }^3=(\xi \lrcorner \;\varphi )^3=\xi \lrcorner \;
[\;(\xi \lrcorner \; \varphi ) \wedge (\xi \lrcorner  \;\varphi )\wedge  \varphi \;]=\xi \lrcorner \; (6 |\xi |^2 \mu )= 6 \mu_{\xi}
\end{equation*}
where $\mu_{\xi}=\mu |_{V_{\xi}}$ is the induced orientation form on $V_{\xi}$.
For  $u,v\in V_{\xi}$.
\begin{eqnarray*}
\omega_{\xi}(J_{\xi}(u), v )&=&\omega_{\xi} ( u\times \xi ,v)=
\langle \psi(u \times \xi,v),\xi\rangle =\varphi (u \times \xi,v, \xi)  
\\
&=&-\varphi (\xi, \;\xi \times u,\;v)= -\langle \;\xi\times (\xi\times u), v \; \rangle \\
&=& -\langle \; - | \xi |^2 u + \langle \xi,u \rangle \xi, v \;\rangle =  | \xi |^2 \langle u,v\rangle - \langle \xi,u \rangle \langle \xi,v\rangle \\
&=& \langle u,v\rangle .
\end{eqnarray*}
implies
$\langle J_{\xi}(u), J_{\xi}(v)\rangle =-\omega_{\xi}(u,J_{\xi}(v))=\langle u,v\rangle$.  By a calculation  of $J_{\xi}$, one checs that the $3$-form
$\Omega_{\xi} $ is a  $(3,0)$ form, furthermore it is non-vanishing because
\begin{eqnarray*} \frac{1}{2i} \;\Omega_{\xi} \wedge \overline{\Omega}_{\xi}=Im\;\Omega_{\xi} \wedge Re\;\Omega_{\xi} &=& (\xi \lrcorner *\varphi)\wedge [\;\xi\lrcorner \;(\xi^{\#}\wedge \varphi)\;]\\
&=&-\xi \lrcorner \; [\;(\xi \lrcorner  *\varphi)\wedge (\xi^{\#}\wedge \varphi)\;]\\
&=&\xi \lrcorner \;[ *(\xi^{\#}\wedge \varphi )\wedge ( \xi^{\#}\wedge \varphi ) \;]\;\;\;\;\;\;\;\\
&=& |\xi^{\#}\wedge \varphi |^2 \; \xi \lrcorner \; vol (M)\\
&=& 4|\xi^{\#}|^2 \; (*\xi^{\#}) =4 \;vol (X_{\xi} ). \;\;\;\;\;\;\;\;
\end{eqnarray*}

\vspace{.1in}

It is easy to see  $*Re \;\Omega_{\xi}=-Im \;
\Omega_{\xi}\wedge \xi^{\#}$ and  $*Im \;\Omega_{\xi}=Re \;
\Omega_{\xi}\wedge \xi^{\#}$. 
\begin{equation*}
*_{3} Re \; \Omega_{\xi}=Im\; \Omega_{\xi}.
\end{equation*}

\vspace{.05in}

Notice that $\omega_{\xi} $ is a symplectic structure on $X_{\xi}$
when  $d \varphi =0$  and $\CL_{\xi}(\varphi)|_{V_{\xi}} =0$,
($\CL_{\xi}$ is the Lie derivative along $\xi$), since $\omega_{\xi} ={\xi} \lrcorner \;\varphi $ and:
$$ d\omega_{\xi}=\CL_{\xi}(\varphi) - \xi  \lrcorner \; d\varphi = \CL_{\xi}(\varphi )$$

\vspace{.1in}

\noindent $J_{\xi}$ is integrable complex structure  if $d^{*}\varphi =0 $ and $\CL_{\xi}(*\varphi)|_{V_{\xi}} =0 $ since 
$$d (Im \Omega_{\xi})=d(\xi \lrcorner *\varphi )=\CL_{\xi }(*\varphi )-\xi \lrcorner \; d(*\varphi )=0$$
Also notice that $d\varphi=0 \implies $ $d
(Re\;\Omega_{\xi})=d(\varphi|_{X_{\xi}})=0$.

\section{HL manifolds and Mirror duality in $G_2$ manifolds}

By \cite{t} any $7$-dimensional Riemanninan manifold admits a non-vanishing orthonormal $2$-frame field $\Lambda=<u, v>$, in particular  $(M, \varphi)$ admits such a field.  $\Lambda $ gives a section of the bundle of oriented $2$-frames $V_{2}(M)\to M$, and hence  gives an associative/coassociative splitting of the tangent bundle $TM={\bf E}\oplus {\bf V}$, where ${\bf E} = {\bf E}_{\Lambda}= <u,v, u\times v>$ and ${\bf V}={\bf V}_{\Lambda}={\bf E}^{\perp}$. When there is no danger of confusion we will denote the $2$-frame fields  and the $2$-planes fields which they induce by the same symbol $\Lambda$. Also,  any unit section $\xi$ of ${\bf E}\to M$ induces a complex structure $J_{\xi}$  on the bundle ${\bf V}\to M$ by the cross product $J_{\xi} (u)=u \times \xi $.

\vspace{.1in}

In \cite{as2} any two almost Calabi-Yau's $X_{\xi}$ and $X_{\xi'}$ inside $(M,\varphi )$ were called {\it dual} if  the defining vector fields $\xi$  and $\xi' $  are chosen from  ${\bf V}$ and ${\bf E}$, respectively. Here we make this correspondence more precise, in particular showing how to choose $\xi$ and $\xi'$  in a more canonical way.

\begin{Def} A $3$-dimensional submanifold  $Y^{3}\subset (M,\varphi)$ is called  {\it Harvey-Lawson manifold} (HL in short) if $\varphi|_{Y}=0$. 
\end{Def}






\vspace{.05in}


\begin{Def}  (\cite {as1}) Call any orthonormal $3$-frame field 
$\Gamma=<u,v,w>$ on $(M,\varphi )$, a {\it $G_2$-frame field} if  $\varphi(u,v,w)= <u\times v, w>=0$, equivalently $w$ is a unit section of ${ \bf V}_{\Lambda} \to X $, with $\Lambda =<u,v>$ (see (1)). 
\end{Def} 


Now pick a nonvanishing $2$-frame field $\Lambda=<u,v>$ on $M$ and 
 let $TM= {\bf E}\oplus {\bf V}$  be the induced splitting with ${\bf E}=<u,v, u\times v>$.  Let $w$ be a unit section of the bundle ${\bf V} \to M$. Such a section $w$ may not exist on whole $M$, but by obstruction theory it exists on a tubular neighborhood $U$ of the $3$-skeleton $M^{(3)}$ of $M$ (which is a complement of some $3$-complex $Z\subset M$). So $\varphi(u,v, w)=0$, and hence $\Gamma=<u,v,w>$ is a $G_{2}$ frame field.
Next consider the non-vanishing vector fields: 

\begin{itemize}
\item $R=\chi(u,v,w)=-u\times(v\times w)$
\item $R'=\frac{1}{\sqrt{3}}(u\times v +v\times w +w\times u)$
\item $R''=\frac{1}{\sqrt{3}}(u +v +w)$
\end{itemize}
If the $6$-plane fields  $R^{\perp}$, $R'^{\perp}$, and  $R''^{\perp}$, 
are integrable we get almost Calabi-Yau manifolds 
$(X_{R}, w_{R}, \Omega_{R}, J_{R})$, 
$(X_{R'}, w_{R'}, \Omega_{R'}, J_{R'})$, and
$(X_{R''}, w_{R''}, \Omega_{R''}, J_{R''})$.
Let us use the convention that $a,b,c$ are real numbers, and  $[u_1,..u_n]$ is the distribution generated by the vectors $u_{1},..,u_{n}$.

 \vspace{.1in}

\begin{Lem} By definitions, the following hold
\begin{itemize}
\item[(a)] $Y:=[u,v,w]=[au+bv+cw\;|\;  a+b+c=0]\oplus[R'']$
\item[(b)]  $\V=[u,v,w, R]$,  is a coassociative $4$-plane field. 
\item[(c)] $ \E: = [u\times v, v\times w, w\times u ]$ is an associative $3$-plane field.
\item[(d)] $\E \perp \V$
\end{itemize}

\end{Lem}

\begin{Thm} For $(a,b,c)\in\R^3$ with  $a+b+c=0$, then 
\begin{itemize}
\item[(a)] $TX_{R}= [au+bv+cw, R'', R', a(v\times w)+b(w\times u)+c(u\times v)]$
\item[] $J_{R}(au+bv+cw) =-a(v\times w)-b(w\times u)-c(u\times v)$
\item[] $J_{R}(R'')=-R'$\\
\item[(b)] $TX_{R'}= [au+bv+cw, R'', R, a(v\times w)+b(w\times u)+c(u\times v)]$
\item[] $J_{R'}(au+bv+cw) =-((b-c)u+(c-a)v+(a-b)w)/\sqrt{3}$
\item[] $J_{R'}(a(v\times w)+b(w\times u) +c(u\times v)) =
 \\
 ((b-c) (v\times w)+(c-a)(w\times u)+ (a-b) (u\times v)   )/\sqrt{3}$
 \item[] $J_{R'}(R'')=R$\\
 \item[(c)] $TX_{R''}= [au+bv+cw, R, R', (a(v\times w)+b(w\times u)+c(u\times v)]$
\item[] $J_{R''}(au+bv+cw) =  \\
((b-a)(u\times v)+ (c-b)(v\times w)+(a-c)(w\times u))/\sqrt{3}$\\
 $J_{R''}(R)=R'$\\
 \item[(d)] $\{u,v,w, R, u\times v, v\times w, w\times u\}$ is an orthonormal frame field.

\end{itemize}
\end{Thm}

\proof To show (a) by using (4) we calculate: \begin{eqnarray*}
R\times u &=& \chi(u,v,w) \times u= -[u\times (v\times w)]\times u =u \times[u\times (v\times w)] \\
&=& -\chi(u,u,v\times w)-<u,u>(v\times w) +<u,v\times w>u \\
&=& -(v\times w) +\varphi(u,v,w)u\\
\end{eqnarray*}
\vspace{-.4in}
\begin{eqnarray}\mbox{Therefore} \;\; R\times u&=&-(v\times w)
\end{eqnarray}
Similarly,  $R\times v=-(w \times u)$ and $R\times w=- (u\times v)$. Therefore we have $J_{R}(au +b v +cw) =-a(v\times w) -b (w\times u)-c(u \times v)$, and $J_{R}(R'')=-R'$.

\begin{eqnarray*}
\sqrt{3} \;R'\times u &=& (u\times v+ v\times w + w \times u)\times u  \\
&=&-u \times (u\times v) -u\times (v\times w)-u\times (w\times u)\\
&=&<u,u>v -<u,v>u \\ && + \;\chi(u,v,w) +<u,v>w-<u,w>v \\&&+<u,w>u-<u,u>w
\end{eqnarray*}
\begin{eqnarray}
\mbox{Therefore} \;\;  \sqrt{3} \;R'\times u &=& R + (v-w)
\end{eqnarray}
Similarly $\sqrt{3} \;R'\times v = R +(w-u)$, and 
$\sqrt{3} \;R'\times w  = \;R +(u-v)$, which implies the first part of (b), and $J_{R'}(R'')=R$. 

\vspace{.1in}


For the second part of (b) we need the compute the following:
\begin{equation}
\sqrt{3}R'\times [a(v\times w) + b(w \times u ) + c(u \times v)] =
\end{equation}
\begin{equation*}
( u\times v +  v\times w + w\times u) \times [a(v\times w) + b(w \times u ) + c(u \times v)]
\end{equation*}


\vspace{.1in}
For this first by repeatedly using (4) and $\varphi(u,v,w)=0$ we calculate:
\begin{eqnarray*}
(v\times u)\times (w\times v) &=&  -\chi(v \times u, w, v)-<v \times u, w>v +<v\times u,v>w\\
&=& -\chi(v \times u, w, v) = -\chi( w, v, v \times u)\\
&=& w \times (v\times (v \times u))+<w, v>(v\times u)-<w,v\times u>v\\
&=&w \times (v\times (v \times u))\\
&=& w\times (-\chi(v,v,u)-<v,v>u +<v,u>v) =- (w\times u)
\end{eqnarray*}
Then by plugging in (9) gives (b). 
Proof of (c) is similar to (a) \qed

\vspace{.15in}

In particular from the above calculations we get can express $\varphi$ as:

\begin{Cor}

$$\varphi = u^{\#}\wedge v^{\#}\wedge (u^{\#}\times v^{\#}) + v^{\#}\wedge w^{\#} \wedge (v^{\#} \times w^{\#}) + w^{\#} \wedge u^{\#} \wedge (w^{\#} \times u^{\#}) $$
$$\;\;\;\;\;\;\; + \;u^{\#} \wedge R^{\#} \wedge (v^{\#}\times w^{\#}) + (v^{\#}\wedge R^{\#}) \wedge (w^{\#} \times u^{\#}) + w^{\#} \wedge R^{\#} \wedge (u^{\#} \times v^{\#}) $$
$$ -\; (u^{\#}\times v^{\#}) \wedge (v^{\#}\times w^{\#} ) \wedge (w^{\#}\times u^{\#}) $$
\end{Cor}

\vspace{.1in}


 
 
\vspace{.05in}

Recall that in an earlier paper we proved the following facts:

{\Prop \cite{as2} Let  $\{ \alpha, \beta \}$ be orthonormal vector fields
on $(M, \varphi) $. Then on $X_{\alpha}$ the following hold
\vspace{.05in}
\begin{itemize}
\item [(i)]$ Re\;\Omega_{\alpha }=\omega_{\beta} \wedge \beta^{\#} + Re\; \Omega_{\beta} $\\
\item [(ii)]$ Im\;\Omega_{\alpha }=\alpha  \lrcorner \; (\star \omega_{\beta})-(\alpha  \lrcorner \; Im\; \Omega_{\beta} )\wedge \beta^{\#} $\\
\item[(iii)] $\omega_{\alpha }= \alpha  \lrcorner \;Re \; \Omega_{\beta } + (\alpha  \lrcorner \; \omega_{\beta })\wedge \beta^{\#} $
\end{itemize} }

\proof  Since $Re\;\Omega_{\alpha } =\varphi |_{X_{\alpha}}$  (i) follows.
Since $Im\; \Omega_{\alpha}=\alpha  \lrcorner *\varphi$  following gives (ii)
\begin{eqnarray*}
\alpha  \lrcorner \; (\star \omega_{\beta})&=& \alpha  \lrcorner \; [\; \beta  \lrcorner \;*(\beta  \lrcorner \;\varphi )\;] \\
&=&\alpha  \lrcorner \; \beta  \lrcorner \; (\beta^{\#} \wedge * \varphi )\\
&=& \alpha  \lrcorner  *\varphi + \beta^{\#}\wedge (\alpha  \lrcorner \;\beta  \lrcorner  *\varphi) \\
&=&  \alpha  \lrcorner  *\varphi + (\alpha  \lrcorner \; Im\; \Omega_{\beta} )\wedge \beta^{\#}
\end{eqnarray*}

(iii) follows from the following computation
$$\alpha  \lrcorner \; Re\; \Omega_{\beta}= \alpha  \lrcorner \;\beta  \lrcorner \; (\beta ^{\#} \wedge \varphi)= \alpha  \lrcorner \; \varphi +\beta^{\#}\wedge (\alpha  \lrcorner \;\beta  \lrcorner \;\varphi )= \alpha  \lrcorner \; \varphi- (\alpha  \lrcorner \; \omega_{\beta })\wedge \beta^{\#} $$
\qed

\vspace{.1in}

Note that even though the identities of this proposition hold only
after restricting the right hand side to $X_{\alpha}$, all the
individual terms are defined everywhere on $(M,\varphi)$. Also,
from the construction, $X_{\alpha }$ and  $X_{\beta }$ inherit
vector fields $\beta$ and $\alpha$, respectively.

\vspace{.05in}

{\Cor \cite{as2} Let  $\{ \alpha, \beta \}$ be orthonormal vector fields on
$(M,\varphi)$. Then there are  $A_{\alpha \beta} \in \Omega^{3}(M)$, and $W_{\alpha \beta}\in \Omega^{2}(M)$ satisfying
\begin{eqnarray*}
(a)\;\;\;\; \;\;\; \; \;\;\; \varphi |_{X_{\alpha}}= Re\;\Omega_{\alpha} &\mbox {and} & \varphi |_{X_{\beta}}=Re\; \Omega_{\beta}\\
(b)\;\;\;\; \;\; A_{\alpha \beta}|_{X_{\alpha}}= Im\;\Omega_{\alpha} &\mbox {and} &
A_{\alpha \beta}|_{X_{\beta}}= \alpha  \lrcorner \; (\star \omega_{\beta})\\
(c) \hspace{.5in} W_{\alpha \beta}|_{X_{\alpha}}= \omega_{\alpha} &\mbox {and} &
W_{\alpha \beta}|_{X_{\beta}}= \alpha  \lrcorner \;Re \; \Omega_{\beta }
\end{eqnarray*} }

Now we can choose $\alpha$ as $R$ and $\beta$ as $R'$ of the given $HL$ manifold. That concludes that given a $HL$ submanifold of a $G_2$ manifold, it will determine a ``canonical'' mirror pair of Calabi-Yau manifolds (related thru the HL manifold) with the complex and symplectic structures given above.

\end{document}